\documentclass[english]{amsart}
\usepackage[T1]{fontenc}
\usepackage[latin9]{inputenc}
\usepackage{babel}
\usepackage{amsthm}
\usepackage{amssymb}
\usepackage[unicode=true,pdfusetitle,
 bookmarks=true,bookmarksnumbered=false,bookmarksopen=false,
 breaklinks=false,pdfborder={0 0 1},backref=false,colorlinks=false]
 {hyperref}

\makeatletter
\numberwithin{equation}{section}
\theoremstyle{plain}
\newtheorem{thm}{\protect\theoremname}[section]
\theoremstyle{plain}
\newtheorem{lem}{\protect\lemmaname}[section]

\makeatother

\providecommand{\lemmaname}{Lemma}
\providecommand{\theoremname}{Theorem}

\begin{document}
\title[Two Lambert series identities]{Proofs for two Lambert series identities of Gosper }
\author{Bing He}
\address{School of Mathematics and Statistics, Central South University\\
Changsha 410083, Hunan, People's Republic of China}
\email{yuhe001@foxmail.com; yuhelingyun@foxmail.com}
\begin{abstract}
Applying the theory of modular forms and Lambert series manipulations
we establish an Eisenstein series identity. From this formula we confirm
a Lambert series identity conjectured by Gosper. Another Lambert series
identity of Gosper is also confirmed by using Lambert series manipulations.

\end{abstract}

\keywords{Lambert series; Modular form; Eisenstein series}
\subjclass[2000]{11F11, 11M36.}
\maketitle

\section{Introduction}

In \cite[pp. 102--103]{G} Gosper conjectured the following three
identities on Lambert series:

\begin{align}
\sum_{n\geq1}\frac{q^{2n-1}}{(1-q^{2n-1})^{2}}-2\sum_{n\geq1}\frac{q^{4n-2}}{(1-q^{4n-2})^{2}} & =\sum_{n\geq1}\frac{(2n-1)q^{2n-1}}{1-q^{4n-2}},\label{eq:tt-1}\\
\sum_{n\geq1}\frac{nq^{2n}}{1+q^{2n}}\sum_{n\geq1}\frac{(2n-1)q^{2n-1}}{1-q^{4n-2}} & =\sum_{n\geq1}\frac{B_{3}(n)}{3}\frac{q^{2n-1}}{1-q^{4n-2}},\label{eq:tt-2}
\end{align}
and
\begin{equation}
6\sum_{n\geq1}\frac{q^{4n-2}}{(1-q^{2n-1})^{4}}+\sum_{n\geq1}\frac{q^{2n-1}}{(1-q^{2n-1})^{2}}=\sum_{n\geq1}\frac{n^{3}q^{n}}{1-q^{2n}},\label{eq:tt-3}
\end{equation}
where 
\[
B_{3}(n):=n(n-1)(n-1/2).
\]

Both sides of these identities involve only Lambert series, but the
second one is different from the other two formulas because the left
side of \eqref{eq:tt-2} is a product of two Lambert series. It should
be emphasized that there are no identities of this type in Gosper's
list except these three formulas. Actually, the identity \eqref{eq:tt-1}
was confirmed by El Bachraoui in \cite[p. 7]{El1} while the other
two formulas remain open.

In this paper we will confirm the Lambert series identities \eqref{eq:tt-2}
and \eqref{eq:tt-3}.
\begin{thm}
\label{T1.1} The identities \eqref{eq:tt-2} and \eqref{eq:tt-3}
hold for $|q|<1.$
\end{thm}
In Section \ref{sec:2}, we employ the theory of modular forms to
establish an Eisenstein series identity, which is crucial in the derivation
of \eqref{eq:tt-2}. Section \ref{sec:3} is devoted to our proof
of Theorem \ref{T1.1}.

\section{\label{sec:2} An auxiliary result}

We first recall some notations from the theory of modular forms.

The modular group $\mathrm{SL}_{2}(\mathbb{Z})$ is defined by 
\[
\mathrm{SL}_{2}(\mathbb{Z}):=\left\{ \left(\begin{matrix}a & b\\
c & d
\end{matrix}\right):a,b,c,d\in\mathbb{Z},\,ad-bc=1\right\} .
\]
If $N>1$ is an integer, the congruence subgroup $\Gamma_{0}(N)$
is defined by
\[
\Gamma_{0}(N)=\left\{ \left(\begin{matrix}a & b\\
c & d
\end{matrix}\right)\in\mathrm{SL}_{2}(\mathbb{Z}):c\equiv0\,(\bmod N)\right\} .
\]
We denote by $M_{k}(\Gamma)$ the space of modular forms of weight
$k$ for $\Gamma$ if $\Gamma$ is a subgroup of $\mathrm{SL}_{2}(\mathbb{Z})$
with finite index.

In order to prove \eqref{eq:tt-2} we need an auxiliary result.
\begin{lem}
For $\mathrm{Im}\:\tau>0,$ we have 
\begin{equation}
\begin{aligned} & E_{4}(\tau)-9E_{4}(2\tau)+8E_{4}(4\tau)\\
 & =10(2E_{2}(4\tau)-E_{2}(2\tau))(3E_{2}(2\tau)-E_{2}(\tau)-2E_{2}(4\tau)),
\end{aligned}
\label{eq:1-4}
\end{equation}
where $E_{2}(\tau)$ and $E_{4}(\tau)$ are two Eisenstein series
given by
\begin{align*}
E_{2}(\tau) & :=1-24\sum_{n\geq1}\sigma_{1}(n)e^{2\pi ni\tau},\\
E_{4}(\tau) & :=1+240\sum_{n\geq1}\sigma_{3}(n)e^{2\pi ni\tau},
\end{align*}
with 
\[
\sigma_{\alpha}(n):=\sum_{d|n}d^{\alpha}.
\]
\end{lem}
\begin{proof}
It follows from \cite[Exercise 1.2.8(e)]{DS} that
\[
E_{2}(\tau)-NE_{2}(N\tau)\in M_{2}(\Gamma_{0}(N)).
\]
In particular,
\[
E_{2}(\tau)-2E_{2}(2\tau)\in M_{2}(\Gamma_{0}(2)).
\]
Then
\[
2E_{2}(4\tau)-E_{2}(2\tau)\in M_{2}(\Gamma_{0}(4))
\]
and
\[
3E_{2}(2\tau)-E_{2}(\tau)-2E_{2}(4\tau)=2E_{2}(2\tau)-E_{2}(\tau)+E_{2}(2\tau)-2E_{2}(4\tau)\in M_{2}(\Gamma_{0}(4)).
\]
Therefore, the right side of \eqref{eq:1-4} belongs to $M_{4}(\Gamma_{0}(4)).$
Since $E_{4}(\tau)\in M_{4}(\mathrm{SL}_{2}(\mathbb{Z})),$ we see
that $E_{4}(2\tau)\in M_{4}(\Gamma_{0}(2))$ and $E_{4}(4\tau)\in M_{4}(\Gamma_{0}(4)).$
From \cite[Theorem 3.5.1]{DS} we know that $\mathrm{dim}(M_{4}(\Gamma_{0}(4)))=3.$
It is clear that $\{E_{4}(\tau),E_{4}(2\tau),E_{4}(4\tau)\}$ is a
basis of $M_{4}(\Gamma_{0}(4)).$ Then there exist three complex numbers
$a,b,c$ such that 
\begin{equation}
\begin{aligned} & 10(2E_{2}(4\tau)-E_{2}(2\tau))(3E_{2}(2\tau)-E_{2}(\tau)-2E_{2}(4\tau))\\
 & =aE_{4}(\tau)+bE_{4}(2\tau)+cE_{4}(4\tau).
\end{aligned}
\label{eq:3-1}
\end{equation}
Set $q=e^{2\pi i\tau}$. Comparing the coefficients of $q^{0},q^{1},q^{2}$
on both sides of \eqref{eq:3-1} we get
\[
a=1,b=-9,c=8.
\]
Then \eqref{eq:1-4} follows readily by substituting these constants
back into \eqref{eq:3-1}.
\end{proof}

\section{\label{sec:3} Proof of Theorem \ref{T1.1}}

We first prove \eqref{eq:tt-2}. Since 
\begin{align*}
\frac{q^{2n}}{1+q^{2n}} & =\frac{q^{2n}}{1-q^{2n}}-\frac{2q^{4n}}{1-q^{4n}},\\
\frac{q^{n}}{1-q^{2n}} & =\frac{q^{n}}{1-q^{n}}-\frac{q^{2n}}{1-q^{2n}}.
\end{align*}
we see that
\begin{equation}
\sum_{n\geq1}\frac{nq^{2n}}{1+q^{2n}}=\sum_{n\geq1}\frac{nq^{2n}}{1-q^{2n}}-2\sum_{n\geq1}\frac{nq^{4n}}{1-q^{4n}},\label{eq:1-1}
\end{equation}
\[
\sum_{n\geq1}\frac{nq^{n}}{1-q^{2n}}=\sum_{n\geq1}\frac{nq^{n}}{1-q^{n}}-\sum_{n\geq1}\frac{nq^{2n}}{1-q^{2n}}
\]
and 
\[
\sum_{n\geq1}\frac{nq^{2n}}{1-q^{4n}}=\sum_{n\geq1}\frac{nq^{2n}}{1-q^{2n}}-\sum_{n\geq1}\frac{nq^{4n}}{1-q^{4n}}.
\]
From the last two identities above we deduce that
\begin{equation}
\begin{aligned}\sum_{n\geq1}\frac{(2n-1)q^{2n-1}}{1-q^{4n-2}} & =\sum_{n\geq1}\frac{nq^{n}}{1-q^{2n}}-\sum_{n\geq1}\frac{2nq^{2n}}{1-q^{4n}}\\
 & =\sum_{n\geq1}\frac{nq^{n}}{1-q^{n}}-3\sum_{n\geq1}\frac{nq^{2n}}{1-q^{2n}}+2\sum_{n\geq1}\frac{nq^{4n}}{1-q^{4n}}.
\end{aligned}
\label{eq:1-2}
\end{equation}
Similarly,
\begin{equation}
\begin{aligned} & \sum_{n\geq1}\frac{B_{3}(n)}{3}\frac{q^{2n-1}}{1-q^{4n-2}}\\
 & =\frac{1}{24}\sum_{n\geq1}(2n-1)((2n-1)^{2}-1)\frac{q^{2n-1}}{1-q^{4n-2}}\\
 & =\frac{1}{24}\sum_{n\geq1}n(n^{2}-1)\frac{q^{n}}{1-q^{2n}}-\frac{1}{12}\sum_{n\geq1}n(4n^{2}-1)\frac{q^{2n}}{1-q^{4n}}\\
 & =\frac{1}{24}\sum_{n\geq1}n(n^{2}-1)\frac{q^{n}}{1-q^{n}}-\frac{1}{24}\sum_{n\geq1}n(n^{2}-1)\frac{q^{2n}}{1-q^{2n}}\\
 & \;-\frac{1}{12}\sum_{n\geq1}n(4n^{2}-1)\frac{q^{2n}}{1-q^{2n}}+\frac{1}{12}\sum_{n\geq1}n(4n^{2}-1)\frac{q^{4n}}{1-q^{4n}}\\
 & =\frac{1}{24}\sum_{n\geq1}n(n^{2}-1)\frac{q^{n}}{1-q^{n}}-\frac{1}{8}\sum_{n\geq1}(3n^{3}-n)\frac{q^{2n}}{1-q^{2n}}\\
 & \;\;+\frac{1}{12}\sum_{n\geq1}n(4n^{2}-1)\frac{q^{4n}}{1-q^{4n}}.
\end{aligned}
\label{eq:1-3}
\end{equation}

Let $q=e^{2\pi i\tau}$ with $\mathrm{Im}\:\tau>0.$ We rewrite the
Eisenstein series $E_{2}(\tau)$ and $E_{4}(\tau)$ as 
\begin{align*}
E_{2}(\tau) & =1-24\sum_{n\geq1}\frac{nq^{n}}{1-q^{n}},\\
E_{4}(\tau) & =1+240\sum_{n\geq1}\frac{n^{3}q^{n}}{1-q^{n}}.
\end{align*}
Then 
\begin{align*}
\sum_{n\geq1}\frac{nq^{n}}{1-q^{n}} & =\frac{1}{24}(1-E_{2}(\tau)),\\
\sum_{n\geq1}\frac{n^{3}q^{n}}{1-q^{n}} & =\frac{1}{240}(E_{4}(\tau)-1).
\end{align*}
Substituting these identities into \eqref{eq:1-1}, \eqref{eq:1-2}
and \eqref{eq:1-3} and then simplifying we get
\begin{align}
\sum_{n\geq1}\frac{nq^{2n}}{1+q^{2n}} & =\frac{1}{24}(2E_{2}(4\tau)-E_{2}(2\tau))-\frac{1}{24},\label{eq:1-11}\\
\sum_{n\geq1}\frac{(2n-1)q^{2n-1}}{1-q^{4n-2}} & =\frac{1}{24}(3E_{2}(2\tau)-E_{2}(\tau)-2E_{2}(4\tau)),\label{eq:1-12}\\
\sum_{n\geq1}\frac{B_{3}(n)}{3}\frac{q^{2n-1}}{1-q^{4n-2}} & =\frac{1}{5760}(E_{4}(\tau)-9E_{4}(2\tau)+8E_{4}(4\tau))\label{eq:1-13}\\
 & \;+\frac{1}{24^{2}}(E_{2}(\tau)-3E_{2}(2\tau)+2E_{2}(4\tau)).\nonumber 
\end{align}
Then \eqref{eq:tt-2} follows easily by adding $\frac{1}{24^{2}}(E_{2}(\tau)-3E_{2}(2\tau)+2E_{2}(4\tau))$
into both sides of \eqref{eq:1-4} and then substituting \eqref{eq:1-11},
\eqref{eq:1-12} and \eqref{eq:1-13} into the resulting identity.

We now show \eqref{eq:tt-3}. It follows easily from the formula 
\[
\frac{1}{1-x}=\sum_{n\geq0}x^{n},\;|x|<1
\]
that 
\begin{equation}
\frac{x}{(1-x)^{2}}=\sum_{n\geq1}nx^{n},\;|x|<1.\label{eq:1-6}
\end{equation}
Differentiating \eqref{eq:1-6} with respect to $x$ and then multiplying
the resulting identity by $x$ we get
\[
\sum_{n\geq1}n^{2}x^{n}=\frac{x(1+x)}{(1-x)^{3}},\;|x|<1.
\]
Differentiating this equation with respect to $x$ again and then
multiplying the resulting identity by $x$ we arrive at
\[
\sum_{n\geq1}n^{3}x^{n}=\frac{x^{3}+4x^{2}+x}{(1-x)^{4}},\;|x|<1.
\]
Then
\[
\sum_{m\geq1}\frac{q^{3m}+4q^{2m}+q^{m}}{(1-q^{m})^{4}}=\sum_{m,n\geq1}n^{3}q^{mn}=\sum_{n\geq1}\frac{n^{3}q^{n}}{1-q^{n}}
\]
and so
\[
\begin{aligned}6\sum_{m\geq1}\frac{q^{2m}}{(1-q^{m})^{4}} & =\sum_{m\geq1}\frac{q^{3m}+4q^{2m}+q^{m}}{(1-q^{m})^{4}}-\sum_{m\geq1}\frac{q^{m}}{(1-q^{m})^{2}}\\
 & =\sum_{n\geq1}\frac{n^{3}q^{n}}{1-q^{n}}-\sum_{m\geq1}\frac{q^{m}}{(1-q^{m})^{2}}.
\end{aligned}
\]
Similarly,
\[
6\sum_{m\geq1}\frac{q^{4m}}{(1-q^{2m})^{4}}=\sum_{n\geq1}\frac{n^{3}q^{2n}}{1-q^{2n}}-\sum_{m\geq1}\frac{q^{2m}}{(1-q^{2m})^{2}}.
\]
We combine these two identities to give 
\[
\begin{aligned}6\sum_{n\geq1}\frac{q^{4n-2}}{(1-q^{2n-1})^{4}} & =6\sum_{n\geq1}\frac{q^{2n}}{(1-q^{n})^{4}}-6\sum_{n\geq1}\frac{q^{4n}}{(1-q^{2n})^{4}}\\
 & =\sum_{n\geq1}\frac{n^{3}q^{n}}{1-q^{2n}}-\sum_{m\geq1}\frac{q^{2m-1}}{(1-q^{2m-1})^{2}}.
\end{aligned}
\]
From this we deduce \eqref{eq:tt-3}. This completes the proof of
Theorem \ref{T1.1}.

\section*{Acknowledgement}

 This work was partially supported by the National Natural Science
Foundation of China (Grant No. 11801451) and the Natural Science Foundation
of Hunan Province (Grant No. 2020JJ5682).


\begin{thebibliography}{1}
\bibitem{DS}F. Diamond and J. Shurman, A first course in modular
forms, Graduate texts in mathematics 228. Springer-Verlag, New York,
2005.

\bibitem{El1}M. El Bachraoui, On series identities of Gosper and
integrals of Ramanujan theta function $\psi(q)$, Proc. Amer. Math.
Soc. 147(10)(2019), 4451\textendash 4464.

\bibitem{G}R.W. Gosper, Experiments and discoveries in $q$-trigonometry,
in: F.G. Garvan, M.E.H. Ismail (Eds.), Symbolic Computation, Number
Theory, Special Functions, Physics and Combinatorics, Kluwer, Dordrecht,
Netherlands, 2001, pp.79\textendash 105.
\end{thebibliography}
\end{document}